\documentclass[a4paper, reqno]{amsart}

\usepackage[british]{babel}
\usepackage{amsthm}
\usepackage{amsmath}
\usepackage[all]{xy}

\theoremstyle{plain}
\newtheorem{theorem}[subsection]{Theorem}
\newtheorem{lemma}[subsection]{Lemma}
\newtheorem{proposition}[subsection]{Proposition}

\theoremstyle{definition}

\newtheorem{remark}[subsection]{Remark}

\newcommand{\defn}{\textbf}
\newcommand{\im}{\mathsf{Im\,}}
\renewcommand{\ker}{\mathsf{Ker\,}}
\newcommand{\A}{\mathcal{A}}

\newcommand{\dom}{\mathrm{dom}}
\newcommand{\cod}{\mathrm{cod}}

\newbox\pullbackbox
\setbox\pullbackbox=\hbox{\xy 0;<1mm,0mm>: \POS(4,0)\ar@{-} (0,0) \ar@{-} (4,4)
\endxy}
\def\pullback{\copy\pullbackbox}

\newdir{>>}{{}*!/3.5pt/:(1,-.2)@^{>}*!/3.5pt/:(1,+.2)@_{>}*!/7pt/:(1,-.2)@^{>}*!/7pt/:(1,+.2)@_{>}}
\newdir{ >}{{}*!/-8pt/@{>}}
\newdir{>}{{}*:(1,-.2)@^{>}*:(1,+.2)@_{>}}
\newdir{<}{{}*:(1,+.2)@^{<}*:(1,-.2)@_{<}}

\begin{document}

\title{A note on the ``Smith is Huq'' condition}

\author{Nelson Martins-Ferreira}
\address{Departamento de Matem\'atica, Escola Superior de Tecnologia e Gest\~ao, Instituto Poli\-t\'ecnico de Leiria}
\address{Centre for Rapid and Sustainable Product Development, Portugal}
\email{nelsonmf@estg.ipleiria.pt}

\author{Tim Van~der~Linden}
\address{Centro de Matem\'atica da Universidade de Coimbra, 3001-454 Coimbra, Portugal}
\email{tvdlinde@vub.ac.be}

\thanks{The first author was supported by IPLeiria/ESTG-CDRSP and Funda\c c\~ao para a Ci\^encia e a Tecnologia (grant number SFRH/BPD/4321/2008). The second author was supported by Centro de Matem\'atica da Universidade de Coimbra and by Funda\c c\~ao para a Ci\^encia e a Tecnologia (grant number SFRH/BPD/38797/2007). He wishes to thank the Instituto Poli\-t\'ecnico for its kind hospitality during his stay in Leiria.}

\begin{abstract}
We show that two known conditions which arose naturally in commutator theory and in the theory of internal crossed modules coincide: \emph{every star-multiplicative graph is multiplicative} if and only if \emph{every two effective equivalence relations commute as soon as their normalisations do}. This answers a question asked by George Janelidze.
\end{abstract}



\maketitle

\section*{Introduction}\label{Section-Introduction}
The purpose of this work is to prove that for a semi-abelian category, the following conditions are equivalent:
\begin{enumerate}
\item[(SM)] every star-multiplicative graph is an internal groupoid;
\item[(SH)] two equivalence relations commute if and only if their normalisations commute.
\end{enumerate}
The first condition comes from the study of internal crossed modules. In a semi-abelian category $\A$, the internal crossed modules introduced by Janelidze~\cite{Janelidze} form a category which is equivalent to the category of internal groupoids in $\A$. To define a crossed module of groups, however, less structure is needed: a reflexive graph equipped with a star-multiplication already determines a crossed module. Nevertheless, there exist examples of semi-abelian categories where this is not true. Thus the question arose under which conditions on~$\A$ the star-multiplicative graphs in $\A$ are internal groupoids.

The second condition was first considered by Bourn and Gran in~\cite{BG}. On one hand, there is the commutator of internal (effective) equivalence relations which was introduced by Smith~\cite{Smith} in the context of Mal'tsev varieties and made categorical by Pedicchio~\cite{Pedicchio}. On the other hand, in the article~\cite{Huq}, Huq introduced a commutator for normal subobjects in a context which is roughly equivalent to that of semi-abelian categories. This definition was further studied by several authors, see e.g., \cite{BG} and \cite{Borceux-Bourn}. Since, in any semi-abelian category, there is a bijective correspondence between the normal subobjects of an object and the effective equivalence relations on it, it is natural to ask how the two concepts of commutator correspond to each other. The answer is that commuting equivalence relations induce commuting normal subobjects~\cite[Proposition~3.2]{BG}, but in general, the concepts are not equivalent---not even in a variety of $\Omega$-groups, as the counterexample of digroups shows~\cite{Borceux-Bourn}. On the other hand, it was shown in~\cite{Gran-VdL} that an equivalence relation $R$ on an object $A$ commutes with the largest equivalence relation~$\nabla_{A}$ as soon as the normalisation $k$ of $R$ is \emph{Huq-central}, i.e., as soon as $k$ commutes with the normalisation $1_{A}$ of $\nabla_{A}$. In fact, a result obtained by Gran says that any two equivalence relations of which the normalisations commute and are jointly strongly epic, commute; see~\cite{EverVdLRCT}. Finally, in a category which is, for instance, pointed and strongly protomodular, any two equivalence relations commute if and only if their normalisations commute~\cite{BG}.

We shall prove that (SH) and (SM) are equivalent conditions. We do this in two steps: in the first section we work towards Theorem~\ref{Theorem-SH-is-KCRG} which essentially states that Condition~(SH) may be restricted to a special class of effective equivalence relations: those pairs of effective equivalence relations which are the kernel pairs of the domain and codomain morphisms of a reflexive graph. Under this latter condition Mantovani and Metere studied the relation between Peiffer graphs and groupoids~\cite[Theorem~6.1]{MM}. We follow their intuition in Section~\ref{Section-Star-Multiplication}, where we prove that a reflexive graph carries a star-multiplication if and only if it is a Peiffer graph if and only if the kernels of its domain and codomain morphisms commute (Proposition~\ref{Structures-Coincide}). This is enough to obtain our main result, Theorem~\ref{Main-Theorem}, which states that (SM) is equivalent to (SH).

\section{The ``Smith is Huq'' condition}\label{Section-Smith-is-Huq}
We show that for a pointed protomodular category, the following two conditions are equivalent:
\begin{enumerate}
\item[(SH)] two effective equivalence relations commute as soon as their normalisations do;
\item[(SH')] every reflexive graph of which the kernels of the domain and the codomain morphisms commute is a groupoid.
\end{enumerate}
Condition (SH) is the \emph{Smith is Huq} condition in the title of this section; condition (SH') is well-known to hold, for instance, in the case of groups: recall the analysis of crossed modules given in the final chapter of Mac\,Lane's~\cite{MacLane}.

\subsection{The context} 
In this section we shall work in \emph{pointed protomodular} categories. A category is \defn{pointed} when it has a \defn{zero object}, i.e., an initial object that is also terminal. A pointed category is \defn{Bourn protomodular}~\cite{Bourn1991} when it is finitely complete and the \defn{Split Short Five Lemma} holds: given a commutative diagram
\begin{equation*}\label{SSFL}
\vcenter{\xymatrix@!0@=5em{K[f] \ar[r]^-{\ker f} \ar[d]_-k & A \ar@<-.5ex>[r]_-{f} \ar[d]^-a & B \ar[d]^-b \ar@<-.5ex>[l]_-{s} \\
K[f'] \ar[r]_-{\ker f'} & A' \ar@<-.5ex>[r]_-{f'} & B', \ar@<-.5ex>[l]_-{s'}}}
\end{equation*}
where $bf=f'a$, $s'b=as$, $f s=1_{B}$ and $f' s'=1_{B'}$, the morphisms $k$ and $b$ being isomorphisms implies that $a$ is an isomorphism. (Note that $s'$ is equal to  $asb^{-1}$, so we could avoid mentioning this morphism and the conditions on it.) 

\begin{lemma}\label{Lemma-Pullback-Iso}
Given a commutative diagram
\begin{equation*}
\vcenter{\xymatrix@!0@=5em{K[f] \ar[r]^-{\ker f} \ar[d]_-k & A \ar@<-.5ex>[r]_-{f} \ar[d]^-a & B \ar[d]^-b \ar@<-.5ex>[l]_-{s} \\
K[f'] \ar[r]_-{\ker f'} & A' \ar[r]_-{f'} & B',}}
\end{equation*}
such that $f s=1_{B}$, the morphism $k$ is an isomorphism if and only if the right hand side commutative square $bf=f'a$ is a pullback.\qed
\end{lemma}

Given a split epimorphism and its kernel as in
\[
\xymatrix@!0@=5em{K \ar[r]^-{k} & A \ar@<-.5ex>[r]_-{f} & B \ar@<-.5ex>[l]_-{s}}
\]
the morphism $k$ and the section $s$ are jointly strongly epic; hence $k$ and $s$ are jointly epic~\cite[Lemma~3.1.22]{Borceux-Bourn}, \cite[Lemma~2.2]{Bourn-Gran-CategoricalFoundations}. For instance, such are the product inclusions $\langle 1_{X},0\rangle\colon{X\to X\times Y}$ and $\langle 0,1_{Y}\rangle\colon{Y\to X\times Y}$.

\subsection{Commuting normal monomorphisms} 
A coterminal pair of morphisms
\begin{equation*}\label{Cospan}
\vcenter{\xymatrix@!0@=5em{X \ar[r]^-{k} & A & Y \ar[l]_-{l}}}
\end{equation*}
\defn{commutes (in the sense of Huq)}~\cite{BG,Huq} when there is a (necessarily unique) morphism $\varphi$ such that the diagram
\[
\xymatrix@!0@=3em{ & X \ar[ld]_{\langle 1_{X},0\rangle} \ar[rd]^-{k} \\
X\times Y \ar@{.>}[rr]|-{\varphi} && A\\
& Y \ar[lu]^{\langle 0,1_{Y}\rangle} \ar[ru]_-{l}}
\]
is commutative. 

We shall only consider the case where $k$ and $l$ are normal monomorphisms (i.e., kernels). We are particularly interested in the situation where they are the kernels of the domain and codomain morphisms of a reflexive graph $C=(C_{1},C_{0},d,c,e)$:
\begin{equation*}\label{RG}
\xymatrix@!0@=5em{C_{1} \ar@<1.5ex>[r]^-{d} \ar@<-1.5ex>[r]_-{c} & C_{0}, \ar[l]|-{e}}
\qquad
d e = c e = 1_{C_{0}}
\end{equation*}
and $k=\ker d\colon{X\to C_{1}}$, $l=\ker c\colon{Y\to C_{1}}$. Using Lemma~\ref{Lemma-Pullback-Iso} we may show that when the kernels $k$ and $l$ of the morphisms $d$ and $c$ in a reflexive graph $C=(C_{1},C_{0},d,c,e)$ commute, their domains are isomorphic.

\begin{lemma}\label{Lemma-Kernels-Commute}
Let $k$ and $l$ be induced by a reflexive graph $C$ as above. If $k$ and $l$ commute then the following commutative squares are pullbacks.
\[
\xymatrix@!0@=5em{X\times Y \ar[r]^-{\pi_{X}} \ar[d]_-{\varphi} & X \ar[d]^{h=c k}\\
C_{1} \ar[r]_-{c} & C_{0}}
\qquad\qquad
\xymatrix@!0@=5em{X\times Y \ar[r]^-{\pi_{Y}} \ar[d]_-{\varphi} & Y \ar[d]^{dl}\\
C_{1} \ar[r]_-{d} & C_{0}}
\]
This makes $X$ and $Y$ isomorphic in a strong sense: there exist morphisms $i\colon{X\to Y}$ and $j\colon{Y\to X}$ such that
\[
j i=1_{X},\qquad i j=1_{Y},\qquad c k  j=d l\qquad\text{and}\qquad c k=d l i.
\]
\end{lemma}
\begin{proof}
The left hand side diagram commutes because $\langle1_{X},0\rangle$ and $\langle0,1_{Y}\rangle$ are jointly epimorphic and moreover $c\varphi\langle1_{X},0\rangle=ck=ck \pi_{X} \langle1_{X},0\rangle$ and
\[
c\varphi\langle 0,1_{Y}\rangle=cl=0=ck \pi_{X} \langle0,1_{Y}\rangle.
\]
It is a pullback by Lemma~\ref{Lemma-Pullback-Iso} since the induced morphism between the kernels of $\pi_{X}$ and $c$ is $1_{Y}$. Similarly the right hand side square is a pullback.

The morphism $i\colon{X\to Y}$ is obtained through the universal property of the first pullback as follows. The equality $c eck=ck=h 1_{X}$ gives rise to a morphism $\iota\colon{X\to X\times Y}$ such that $\varphi\iota=eck$ and $\pi_{X} \iota=1_{X}$; considering $X\times Y$ as a product now, this $\iota$ is a pair $\langle1_{X},i\rangle\colon {X\to X\times Y}$. Clearly,
\[
dli=dl\pi_{Y}\langle1_{X},i\rangle=d\varphi\langle1_{X},i\rangle=deck=ck.
\]
Using the second pullback one obtains a morphism $j\colon{Y\to X}$ satisfying $\varphi \langle j,1_{Y}\rangle=edl$, so that $c k  j=d l$. 

Now we only have to prove that $i$ and $j$ are mutually inverse. This again follows from the universal properties of the pullbacks. Indeed, the morphisms $\langle j,ij\rangle\colon {Y\to X\times Y}$ and $\langle j,1_{Y}\rangle\colon {Y\to X\times Y}$ are both universally induced by the equality $cedl=ckj=hj$, hence they are equal. Likewise, $\langle 1_{X},i\rangle$ is equal to $\langle ji,i\rangle$ so that $ji=1_{X}$.
\end{proof}

This result may be interpreted as follows: the two (\emph{a priori} non-equivalent) ways a reflexive graph can be normal\-ised---mapping $C$ to either $ck\colon{X\to C_{0}}$ or $dl\colon{Y\to C_{0}}$---induce naturally isomorphic functors from the category of reflexive graphs with commuting kernels to the category of objects over $C_{0}$.

One usually views the elements of $C_{1}$ as arrows between the elements of $C_{0}$, so that the morphism $\varphi\colon{X\times Y\to C_{1}}$ is nothing but a partial composition on $C_{1}$ which sends a pair of arrows
\[
\xymatrix@1{{} \cdot  & 0 \ar[l]_-{\alpha} & \ar[l]_-{\beta} \cdot}
\]
to its composite $\varphi(\alpha,\beta)$. The central question studied in this paper is under which conditions such a partial composition extends to a composition on the entire graph. To answer it, we shall need the concept of commuting effective equivalence relations and its connection with commuting normal monomorphisms.

\subsection{Commuting effective equivalence relations}
Consider a pair of equivalence relations $(R,S)$ on a common object $A$
\begin{equation*}\label{Category-RG}
\xymatrix@!0@=5em{R \ar@<1.5ex>[r]^-{r_{0}} \ar@<-1.5ex>[r]_-{r_{1}} & A \ar[l]|-{\Delta_{R}} \ar[r]|-{\Delta_{S}} & S, \ar@<1.5ex>[l]^-{s_{0}} \ar@<-1.5ex>[l]_-{s_{1}}}
\end{equation*}
and consider the induced pullback of $r_{1}$ and $s_{0}$.
\begin{equation}\label{Pullback-RS}
\vcenter{\xymatrix@!0@=5em{R\times_{A}S \ar@{}[rd]|<<{\pullback} \ar[r]^-{\pi_{S}} \ar[d]_-{\pi_{R}} & S \ar[d]^-{s_{0}} \\
R \ar[r]_-{r_{1}} & A}}
\end{equation}
The pair $(R,S)$ \defn{commutes (in the sense of Smith)}~\cite{Smith,Pedicchio,BG} when there is a (necessarily unique) morphism $\theta$ such that the diagram
\[
\xymatrix@!0@=3em{ & R \ar[ld]_{\langle 1_{R},\Delta_{S} r_{1}\rangle} \ar[rd]^-{r_{0}} \\
R\times_{A}S \ar@{.>}[rr]|-{\theta} && A\\
& S \ar[lu]^{\langle\Delta_{R} s_{0},1_{S}\rangle} \ar[ru]_-{s_{1}}}
\]
is commutative.

We shall only consider the case where $R$ and $S$ are effective equivalence relations (i.e., kernel pairs). It is well-known that when for a span
\begin{equation}\label{span}
\vcenter{\xymatrix@!0@=3em{ & C_{1} \ar[ld]_-{d} \ar[rd]^-{c} \\
C_{0} && C'_{0},}}
\end{equation}
the kernel pairs $R[d]$ and $R[c]$ commute, this means that $(d,c)$ carries an internal pregroupoid structure~\cite{Janelidze-Pedicchio}; briefly, any zigzag
\[
\xymatrix@1{{} \cdot  & \cdot \ar[l]_-{\alpha} \ar[r]^-{\beta} & \cdot & \cdot \ar[l]_-{\gamma}}
\]
in $C_{1}$ may be composed to a single arrow $\theta(\alpha,\beta,\gamma)$, in such a way that $\theta(\alpha,\beta,\beta)=\alpha$ and $\theta(\beta,\beta,\gamma)=\gamma$. In particular, a reflexive graph $C=(C_{1},C_{0},d,c,e)$ is an internal groupoid if and only if $R[d]$ and $R[c]$ commute: then $\theta(\alpha,\beta,\gamma)=\alpha\circ\beta^{-1}\circ\gamma$.

It is also well-known that when a pair $(R,S)$ of (effective) equivalence relations commutes, then so do their normalisations
\[
\xymatrix@=5em{X=K[r_{0}] \ar[r]^-{k=r_{1} \ker r_{0}} & A & K[s_{0}]=Y: \ar[l]_-{l=s_{1} \ker s_{0}}}
\]
see~\cite[Proposition~3.2]{BG}. In particular, for any internal groupoid $C$ the composition on~$C$ restricts in such a way that the kernels of its domain and codomain morphisms commute. The converse is not true: in general, it is not possible to extend the partial composition on a reflexive graph which is given by its commuting kernels to a composition on the entire graph which makes it into a groupoid. This is explained by the following result (inspired by Lemma~2.1 in~\cite{Johnstone:Herds}), together with the fact that a pair of effective equivalence relations of which the normalisations commute need not commute itself~\cite{Borceux-Bourn}.

\begin{theorem}\label{Theorem-SH-is-KCRG}
For a pointed protomodular category, the following conditions are equivalent:
\begin{enumerate}
\item[(SH)] two effective equivalence relations commute as soon as their normalisations do;
\item[(SH')] every reflexive graph with commuting kernels of the domain and the codomain morphisms is a groupoid.
\end{enumerate}
\end{theorem}
\begin{proof}
It is clear that (SH') is just (SH) in the special case where the effective equivalence relations considered are the kernel pairs of the domain and the codomain morphisms of a reflexive graph. This special case implies the general case. Indeed, let $R=R[d]$ and $S=R[c]$ be the effective equivalence relations induced by a span~\eqref{span} and assume that the normal monomorphisms $k=\ker d$ and $l=\ker c$ commute in the sense of Huq. We have to prove that $R$ and~$S$ commute in the sense of Smith, i.e., the span $(d,c)$ is a pregroupoid. 

If one thinks of the ``elements'' of the object $C_{1}$ as arrows $\xymatrix@1{d(\alpha)\ar[r]^{\alpha} & c(\alpha)}$ then $R$ and~$S$ consist of pairs
\[
\xymatrix@1{{} \cdot & \cdot \ar[r]^-{\alpha} \ar[l]_-{\beta} & \cdot}
\qquad\text{and}\qquad
\xymatrix@1{{} \cdot \ar[r]^-{\gamma} & \cdot & \cdot \ar[l]_-{\delta}},
\]
respectively. Forming the pullback~\eqref{Pullback-RS} of $r_{1}$ and $s_{0}$ we obtain a reflexive graph
\begin{equation}\label{Groupoid}
\xymatrix@=5em{R\times_{C_{1}}S  \ar@<1.5ex>[r]^-{\dom=r_{0}\pi_{R}} \ar@<-1.5ex>[r]_-{\cod=s_{1}\pi_{S}} & C_{1}. \ar[l]|-{\langle\Delta_{R},\Delta_{S}\rangle}}
\end{equation}
An element of $R\times_{C_{1}}S$ is a triple
\[
\xymatrix@1{{} \cdot  & \cdot \ar[l]_-{\alpha} \ar[r]^-{\beta} & \cdot & \cdot \ar[l]_-{\gamma}}
\]
considered as an arrow $\beta$ with domain $\alpha=\dom(\alpha,\beta,\gamma)=r_{0}\pi_{R}(\alpha,\beta,\gamma)$ and codomain $\gamma=\cod(\alpha,\beta,\gamma)=s_{1}\pi_{S}(\alpha,\beta,\gamma)$.
The kernels $\dom$ and $\cod$ commute because so do $k$ and~$l$: the needed morphism
\[
{K[\dom]\times K[\cod]\to R\times_{C_{1}}S}
\]
takes a pair
\[
(\xymatrix@1{{} \cdot  & \cdot \ar[l]_-{0} \ar[r]^-{\beta} & \cdot & \cdot \ar[l]_-{\gamma}}, 
\xymatrix@1{{} \cdot  & \cdot \ar[l]_-{\delta} \ar[r]^-{\epsilon} & \cdot & \cdot \ar[l]_-{0}})
\]
in the product $K[\dom]\times K[\cod]$ and maps it to the element
\[
\xymatrix@1{{} \cdot  & \cdot \ar[l]_-{\delta} \ar[r]^-{\varphi(\beta,\epsilon)} & \cdot & \cdot \ar[l]_-{\gamma}}
\]
of $R\times_{C_{1}}S$. The hypothesis that (SH') holds now implies that this reflexive graph is a groupoid. This, in turn, establishes a pregroupoid structure on the span $(d,c)$: the required morphism $\theta\colon{R\times_{C_{1}}S\to C_{1}}$ is determined by
\[
(\xymatrix@1{{} \cdot  & \cdot \ar[l]_-{\gamma} \ar[rr]^-{\theta(\alpha,\beta,\gamma)} && \cdot & \cdot \ar[l]_-{\alpha}})
=
(\xymatrix@1{{} \cdot  & \cdot \ar[l]_-{\beta} \ar[r]^-{\alpha} & \cdot & \cdot \ar[l]_-{\alpha}})
\circ
(\xymatrix@1{{} \cdot  & \cdot \ar[l]_-{\gamma} \ar[r]^-{\gamma} & \cdot & \cdot \ar[l]_-{\beta}})
\]
where the composition takes place in the groupoid~\eqref{Groupoid}. Indeed, in this groupoid
\[
(\beta,\beta,\beta)\circ (\gamma,\gamma,\beta)=(\gamma,\gamma,\beta)
\]
so that $\theta(\beta,\beta,\gamma)=\gamma$. Likewise, $\theta(\alpha,\beta,\beta)=\alpha$. 
\end{proof}

Condition (SH) is sometimes called the \defn{Smith is Huq} property. It is known to hold in quite diverse situations: in pointed and \emph{strongly protomodular} categories (by~\cite{BG}; see also~\cite{Borceux-Bourn} and~\cite{Bourn2004}) and in pointed and \emph{action accessible} categories (as explained in~\cite{MM}; see also~\cite{BJ07}). This condition is also weaker than the \emph{reflected admissibility} condition studied in~\cite{NMF2}.

\begin{remark}
As explained to us by Tomas Everaert, the condition (SH) may be replaced by its non-effective version
\begin{enumerate}
\item[(SH'')] two equivalence relations commute as soon as their normalisations do,
\end{enumerate}
using the same proof, even when the category is not Barr exact. Then the kernels should be replaced by normal monomorphisms in the sense of Bourn~\cite{Bourn2000}.
\end{remark}

\section{Star-multiplication}\label{Section-Star-Multiplication}
In this section we show that, in a semi-abelian category, three types of (uniquely determined) structure on a reflexive graph $C=(C_{1},C_{0},d,c,e)$ coincide: a reflexive graph $C$ is \emph{star-multiplicative} if and only if it is \emph{Peiffer} if and only if the kernels of $d$ and $c$ commute (Proposition~\ref{Structures-Coincide}). This allows us to prove Theorem~\ref{Main-Theorem} which states that a semi-abelian category has the \emph{Smith is Huq} property if and only if every star-multiplicative graph is a groupoid.

\subsection{The context}
A category is \defn{semi-abelian}~\cite{Janelidze-Marki-Tholen} when it is pointed, Bourn protomodular and Barr exact with binary coproducts. \defn{Barr exact} means that every internal equivalence relation is \defn{effective} (i.e., it is a kernel pair) and the category is \defn{regular}: finitely complete with pullback-stable regular epimorphisms and coequalisers of effective equivalence relations. A \defn{homological} category is pointed, regular and protomodular~\cite{Borceux-Bourn}.

In a homological category regular epimorphisms (coequalisers), strong epimorphisms and normal epimorphisms (cokernels) coincide, and every morphism $f\colon{A\to B}$ may be factored as a regular epimorphism ${A\to I[f]}$ followed by a monomorphism $\im f\colon{I[f]\to B}$. The monomorphism $\im f$ is the \defn{image} of $f$. A morphism $f$ is \defn{proper} when it has a normal image, i.e., $\im f$ is a normal monomorphism. In a semi-abelian category, the \defn{direct image} $\im (p m)$ of a normal monomorphism~$m$ along a regular epimorphism $p$ is always a normal monomorphism (condition (SA*6) in~\cite{Janelidze-Marki-Tholen}). 

We need the following strengthening of Lemma~\ref{Lemma-Pullback-Iso}; see~\cite{Borceux-Bourn} or~\cite[Proposition~7]{Bourn2001}.

\begin{lemma}\label{Lemma-Pullback-Iso-Regular}
In a homological category, given a commutative diagram
\[
\xymatrix@!0@=5em{K[f] \ar[r]^-{\ker f} \ar[d]_-k & A \ar[r]^-{f} \ar[d]^-a & B \ar[d]^-b \\
K[f'] \ar[r]_-{\ker f'} & A' \ar[r]_-{f'} & B'}
\]
where $f$ is a regular epimorphism, the morphism $k$ is an isomorphism if and only if the right hand side square $bf=f'a$ is a pullback.\qed
\end{lemma}

\subsection{Star-multiplicative graphs}
A reflexive graph $C=(C_{1},C_{0},d,c,e)$ is \defn{star-multiplicative}~\cite{Janelidze} when there is a (necessarily unique) morphism
\[
\varsigma\colon{C_{1}\times_{C_{0}} X\to X}
\]
such that $\varsigma \langle k,0\rangle=1_{X}$ and $\varsigma \langle eck,1_{X}\rangle=1_{X}$. Here the square 
\[
\vcenter{\xymatrix@!0@=5em{C_{1}\times_{C_{0}} X \ar[d]_-{\pi_{0}} \ar[r]^-{\pi_{1}} & X \ar[d]^-{h=ck}\\
C_{1} \ar[r]_-{d} & C_{0}}}
\]
is a pullback. A star-multiplication takes a composable pair of arrows
\[
\xymatrix@1{\cdot  & \cdot \ar[l]_-{\alpha} & 0 \ar[l]_-{\beta}}
\]
and sends it to their composite $\varsigma(\alpha,\beta)$.

\subsection{Peiffer graphs}
A reflexive graph $C=(C_{1},C_{0},d,c,e)$ is \defn{Peiffer} when there is a (necessarily unique) morphism
\[
\omega\colon{X\times X\to C_{1}}
\]
such that $\omega\langle 1_{X},0\rangle=k$ and $\omega\langle 1_{X},1_{X}\rangle =eck$. (This definition is not the original one given in~\cite{MM}, but it is equivalent to it in the present context; see~\cite[Theorem 5.3]{MM}.) The structure $\omega$ sends a composable pair of arrows
\[
\xymatrix@1{{}\cdot & 0 \ar[r]^-{\beta} \ar[l]_-{\alpha} & \cdot}
\]
to the composite $\omega(\alpha,\beta)$---which should be considered as $\alpha\circ\beta^{-1}$.

In~\cite{MM} these two structures are shown to be equivalent; we recall the argument.

\begin{proposition}\label{Proposition-s-w-s}
A reflexive graph $C=(C_{1},C_{0},d,c,e)$ in a pointed protomodular category is star-multiplicative if and only if it is Peiffer.
\end{proposition}
\begin{proof}
Given $\varsigma\colon{C_{1}\times_{C_{0}} X\to X}$ put $\omega=\pi_{0}\langle\varsigma,\pi_{1}\rangle^{-1}$; given $\omega\colon{X\times X \to C_{1}}$ put $\varsigma=\pi_{0}\langle\omega,\pi_{1}\rangle^{-1}$. Notations are as above. The inverse morphisms exist by the Split Short Five Lemma.
\end{proof}

Now we work towards an equivalence with reflexive graphs of which the kernel of the domain morphism commutes with the kernel of the codomain morphism. In Lemma~\ref{Lemma-Peiffer} we need the surrounding category to be semi-abelian.

\begin{lemma}\cite[Theorem 5.3]{MM}\label{Lemma-Peiffer-1}
Any Peiffer graph $C$ induces commutative squares
\[
\vcenter{\xymatrix@!0@=5em{X\times X \ar@{}[rd]|-{\mathtt{(i)}} \ar[d]_-{\omega} \ar[r]^-{\pi_{1}} & X \ar[d]^-{h=ck}\\
C_{1} \ar[r]_-{d} & C_{0}}}
\qquad\text{and}\qquad
\vcenter{\xymatrix@!0@=5em{X\times X \ar@{}[rd]|-{\mathtt{(ii)}} \ar[d]_-{\omega} \ar[r]^-{\pi_{0}} & X \ar[d]^-{h=ck}\\
C_{1} \ar[r]_-{c} & C_{0}.}}
\]
Furthermore, the square~$\mathtt{(i)}$ is a pullback.
\end{lemma}
\begin{proof}
The morphisms $\langle1_{X},0\rangle$ and $\langle1_{X},1_{X}\rangle$ are jointly epic and
\begin{align*}
& d\omega\langle1_{X},0\rangle=dk=0=h\pi_{1}\langle1_{X},0\rangle,\\
& d\omega\langle1_{X},1_{X}\rangle=deck=ck=h\pi_{1}\langle1_{X},1_{X}\rangle,\\
& c\omega\langle1_{X},0\rangle=ck=0=ck\pi_{0}\langle1_{X},0\rangle
\end{align*}
and
\[
c\omega\langle1_{X},1_{X}\rangle=ceck=ck=h\pi_{0}\langle1_{X},1_{X}\rangle
\]
so that the two squares commute. Taking kernels horizontally in~$\mathtt{(i)}$ induces the identity morphism $1_{X}$; hence the square is a pullback by Lemma~\ref{Lemma-Pullback-Iso}.
\end{proof}

\begin{lemma}\label{Lemma-Product}
Let $g\colon{X\times X\to A}$ be a morphism with $g\langle 0,1_{X}\rangle=0$ and write $g_{0}=g\langle 1_{X},0\rangle$. Then $g=g_{0} \pi_{0}$, so that $g \langle 1_{X},1_{X}\rangle=g_{0}$.
\end{lemma}
\begin{proof}
The morphism $g$ is uniquely determined by the equalities $g\langle 0,1_{X}\rangle=0$ and $g\langle1_{X},0\rangle=g_{0}$. Since also $g_{0} \pi_{0}\langle0,1_{X}\rangle=0$ and $g_{0} \pi_{0}\langle1_{X},0\rangle=g_{0}$ we have that $g=g_{0} \pi_{0}$.
\end{proof}

\begin{lemma}\label{Lemma-c-cokernel}
For any Peiffer graph $C$, the morphism $c$ is the cokernel of the composite $\omega\langle0,1_{X}\rangle\colon{X\to C_{1}}$.
\end{lemma}
\begin{proof}
First note that $c\omega\langle0,1_{X}\rangle=0$ by commutativity of the square~$\mathtt{(ii)}$ in Lemma~\ref{Lemma-Peiffer-1}. Consider $f\colon{C_{1}\to A}$ with $f \omega\langle0,1_{X}\rangle=0$; we claim that the morphism $fe\colon {C_{0}\to A}$ satisfies $fec=f$. Indeed, by Lemma~\ref{Lemma-Product} the equalities $f \omega\langle0,1_{X}\rangle=0$ and $f\omega\langle 1_{X},0\rangle=fk$ imply $f\omega\langle1_{X},1_{X}\rangle=fk$, so that $feck=fk$. Since also $fece=fe$ and $k$ and $e$ are jointly epic we may conclude that $fec=f$.
\end{proof}

\begin{lemma}\label{Lemma-Peiffer}
For any Peiffer graph $C$ in a semi-abelian category the induced commutative square~$\mathtt{(ii)}$ from Lemma~\ref{Lemma-Peiffer-1} is a pullback.
\end{lemma}
\begin{proof}
Taking kernels vertically gives rise to the reflexive graph
\[
\vcenter{\xymatrix@=3em{K[\omega] \ar@<1.5ex>[r]^-{\pi'_{0}} \ar@<-1.5ex>[r]_-{\pi'_{1}} & K[h]; \ar[l]|-{\Delta}}}
\]
Since~$\mathtt{(i)}$ is a pullback, the morphism $\pi'_{1}$, and hence also $\pi'_{0}$, is an isomorphism. It follows by Lemma~\ref{Lemma-Pullback-Iso-Regular} that the top square in the vertical regular epi-mono factorisation
\[
\xymatrix@!0@=5em{X \ar@{.>}[r]^-{\langle0,1_{X}\rangle} \ar@{:}[d] & X\times X \ar@{}[rd]|<<{\pullback} \ar@{>>}[d] \ar[r]^-{\pi_{0}} & X \ar@{>>}[d]\\
X \ar@{.>}[r]^-{\ker \overline{c}} \ar@{}[rd]|-{\mathtt{(iii)}} \ar@{.>}[d]_-{i} & I[\omega] \ar@{{ >}->}[d]^-{\im \omega} \ar[r]^-{\overline{c}} & I[h] \ar@{{ >}->}[d]^{\im h}\\
Y \ar@{.>}[r]_-{\ker c}& C_{1} \ar[r]_-{c} & C_{0}}
\]
of~$\mathtt{(ii)}$ is a pullback. Taking kernels to the left induces morphisms as indicated. We have to show that $i$ is an isomorphism.

Being a composite $h=ck$ of a normal monomorphism with a regular epimorphism, the morphism $h$ is proper, i.e., its image $\im h$ is a normal monomorphism. Since the square~$\mathtt{(i)}$ is a pullback, $\omega$ is also proper, so that $\im \omega$ is a normal monomorphism. The morphism $\im h$ being mono implies that the square~$\mathtt{(iii)}$ is a pullback. Since both $\im \omega$ and $\ker c$ are normal monomorphisms, this implies that the diagonal of~$\mathtt{(iii)}$---the morphism $\omega \langle0,1_{X}\rangle$---is also a normal monomorphism. Lemma~\ref{Lemma-c-cokernel} tells us that $c$ is its cokernel, so that $\omega \langle0,1_{X}\rangle$ is the kernel of $c$. This means that $i$ is an isomorphism, and the square~$\mathtt{(ii)}$ is a pullback by Lemma~\ref{Lemma-Pullback-Iso}.
\end{proof}

\begin{proposition}\label{Structures-Coincide}
For a reflexive graph $C=(C_{1},C_{0},d,c,e)$ in a semi-abelian category, the following three conditions are equivalent:
\begin{enumerate}
\item $C$ is star-multiplicative;
\item $C$ is Peiffer;
\item $\ker d$ and $\ker c$ commute.
\end{enumerate}
\end{proposition}
\begin{proof}
Conditions (1) and (2) are equivalent by Proposition~\ref{Proposition-s-w-s}. If $C$ is Peiffer then $\ker d$ and $\ker c$ commute. Indeed, by Lemma~\ref{Lemma-Peiffer} we can put $\varphi=\omega$ since $\omega\langle0,1_{X}\rangle$ is the kernel $l$ of~$c$. Conversely, if Condition (3) holds then by Lemma~\ref{Lemma-Kernels-Commute} we have
\[
\iota=\langle1_{X},i\rangle\colon{X\to X\times Y}
\]
such that $\varphi\iota=eck$. Now $\omega=\varphi(1\times i)\colon{X\times X\to C_{1}}$ is a Peiffer structure on $C$ because
\[
\omega\langle 1_{X},0\rangle=\varphi(1_{X}\times i)\langle1_{X},0\rangle=\varphi\langle1_{X},0\rangle=k
\]
and $\omega\langle1_{X},1_{X}\rangle=\varphi(1_{X}\times i)\langle1_{X},1_{X}\rangle=\varphi\iota=eck$.
\end{proof}

\begin{theorem}\label{Main-Theorem}
For a semi-abelian category, the following conditions are equivalent:
\begin{enumerate}
\item[(SM)] every star-multiplicative graph is multiplicative;
\item[(SH)] two (effective) equivalence relations commute if and only if their normalisations commute.
\end{enumerate}
\end{theorem}
\begin{proof}
We already explained above that one implication of (SH) always holds by~\cite[Proposition~2.7.7]{Borceux-Bourn}. Hence by Theorem~\ref{Theorem-SH-is-KCRG} we may replace the second condition with
\emph{\begin{enumerate}
\item[(SH')] every reflexive graph with commuting kernels of the domain and the codomain morphisms is a groupoid.
\end{enumerate}}
The result now follows from Proposition~\ref{Structures-Coincide} and the fact that in a semi-abelian category, multiplicative graphs (i.e., categories) and groupoids coincide.
\end{proof}

Note that Lemma~\ref{Lemma-Peiffer} is the only place where we use that the underlying category is semi-abelian rather than pointed protomodular. This suggests an extension of the concept of Peiffer graph to pointed protomodular categories, where the pullback property of square~$\mathtt{(ii)}$ in Lemma~\ref{Lemma-Peiffer-1} becomes an axiom. (Or, equivalently, in the homological case, the morphism $\omega\langle 0,1_{X}\rangle$ is demanded to be a normal monomorphism.) The concept of star-multiplicative graph allows a similar modification, where now one asks that the morphism of reflexive graphs
\[
\xymatrix@!0@=5em{{C_{1}\times_{C_{0}}X} \ar@<-1.5ex>[d]_-{\pi_{1}} \ar@<1.5ex>[d]^-{\varsigma} \ar[r]^-{\pi_{0}} & C_{1} \ar@<-1.5ex>[d]_-{d} \ar@<1.5ex>[d]^-{c} \\
X \ar[r]_-{h} \ar[u] & C_{0} \ar[u]}
\]
is not just a discrete cofibration (i.e., the square $h\pi_{1}=d\pi_{0}$ is a pullback) but also a discrete fibration ($h\varsigma=c\pi_{0}$ is a pullback). These definitions extend Theorem~\ref{Main-Theorem} to the pointed protomodular context.

\subsection*{Acknowledgements}
Thanks to Tomas Everaert and Julia Goedecke for interesting suggestions and for their comments on the text.

%\bibliography{tim}
%\bibliographystyle{amsalpha}
%% .bbl

\providecommand{\bysame}{\leavevmode\hbox to3em{\hrulefill}\thinspace}
\providecommand{\MR}{\relax\ifhmode\unskip\space\fi MR }
\providecommand{\MRhref}[2]{%
  \href{http://www.ams.org/mathscinet-getitem?mr=#1}{#2}
}
\providecommand{\href}[2]{#2}

\end{document}